\documentclass[12pt,reqno]{amsart}
\usepackage{amssymb}
\usepackage[all]{xy}

\theoremstyle{plain}
\newtheorem*{thm}{Theorem}
\newtheorem*{lem}{Lemma}
\newtheorem*{prop}{Proposition}
\newtheorem*{cor}{Corollary}

\theoremstyle{remark}

\newtheorem*{notat}{Notations and Conventions}

\newfam\cyrfam
\font\tencyr=wncyr10 at 11pt
\font\sevencyr=wncyr7 at 11pt
\font\fivecyr=wncyr5 at 11pt
\newfam\cyrfam
\textfont\cyrfam=\tencyr
\scriptfont\cyrfam=\sevencyr
\scriptscriptfont\cyrfam=\fivecyr
\def\cyr#1{{\fam\cyrfam\relax#1}}

\def\Z{\mathbb{Z}}
\def\Q{\mathbb{Q}}

\def\S{\mathcal{S}}
\def\A{\mathcal{A}}
\def\T{\mathcal{T}}
\def\P{\mathcal{P}_2}
\def\G{\mathcal{G}}
\def\H{\mathcal{H}}
\def\O{\mathcal{O}}
\def\stab{\operatorname{stab}}
\def\rk{\operatorname{rank}}
\def\Im{\operatorname{Im}}

\def\stab{\operatorname{stab}}
\def\Ker{\operatorname{Ker}}
\def\Res{\operatorname{Res}}
\def\Hom{\operatorname{Hom}}
\def\Ext{\operatorname{Ext}}
\def\End{\operatorname{End}}
\def\Aut{\operatorname{Aut}}
\def\PGLn{\operatorname{PGL}_n(k)}
\def\GL{\operatorname{GL}}
\def\Pic{\operatorname{Pic}}
\def\sym#1{\operatorname{\mathsf{Sym}}^2{#1}}
\def\wed#1{\bigwedge^2{#1}}
\def\ind#1#2{\!\!\uparrow_{#1}^{#2}}
\def\gen#1{\langle{#1}\rangle}
\def\tto{\longrightarrow}
\def\onto{\twoheadrightarrow}

\def\iso{\stackrel{\cong}{\longrightarrow}}
\def\iff{\Longleftrightarrow}

\newdir{ >}{{}*!/-5pt/\dir{>}}

\begin{document}

%Topmatter%%%%%%%%%%%%%%%%%%%%%%%%%%%%%%%%

\title[lattices]{On certain lattices associated with \\ generic division algebras}

\author{Nicole Lemire}
\address{Department of Pure Mathematics\\
        University of Waterloo\\
        Waterloo, ON, Canada, N2L 3G1}
\email{nlemire@math.uwaterloo.ca}
%\thanks{}

\author{Martin Lorenz}
\address{Department of Mathematics\\
        Temple University\\
        Philadelphia, PA 19122-6094}
\email{lorenz@math.temple.edu}
\thanks{Research of the second author supported in part by NSF Grant
DMS-9618521}



\begin{abstract}
Let $\S_n$ denote the symmetric group on $n$ letters. We consider the
$\S_n$-lattice $A_{n-1}=\{(z_1,\ldots,z_n)\in\Z^n\mid\sum_iz_i=0\}$, where
$\S_n$ acts on $\Z^n$ by permuting the coordinates, and its squares
$A_{n-1}^{\otimes 2}$, $\sym{A_{n-1}}$, and 
$\wed{A_{n-1}}$. For odd values of $n$, we show that $A_{n-1}^{\otimes 2}$
is equivalent to $\wed{A_{n-1}}$ in the sense of Colliot-Th\'el\`ene and
Sansuc \cite{CTS}. Consequently, the rationality problem for generic division
algebras amounts to proving stable rationality of the multiplicative invariant
field $k(\wed{A_{n-1}})^{\S_n}$ ($n$ odd). Furthermore, confirming a conjecture
of Le Bruyn \cite{LeB}, we show that $n=2$ and $n=3$ are the only cases where
$A_{n-1}^{\otimes 2}$ is equivalent to a permutation $\S_n$-lattice. In the
course of the proof of this result, we construct subgroups $\H\le\S_n$, for all
$n$ that are not prime, so that the multiplicative invariant algebra 
$k[A_{n-1}]^{\H}$ has a non-trivial Picard group.
\end{abstract}

\maketitle

%end Topmatter%%%%%%%%%%%%%%%%%%%%%%%%%%%%%%%%%%

\section*{Introduction} \label{intro}

This note addresses two points raised in Le Bruyn's excellent survey article
\cite{LeB} on the rationality problem for the center of generic division
algebras. This problem, in a reformulation due to Formanek \cite{Fo}, asks
whether the multiplicative invariant field 
$$
k(A_{n-1}^{\otimes 2})^{\S_n}
$$
that is associated with the action of the symmetric group $\S_n$ on the
tensor square of the root lattice $A_{n-1}$ is (stably) rational. 
Recall that $A_{n-1}=\{(z_1,\ldots,z_n)\in\Z^n\mid\sum_iz_i=0\}$, with
$\S_n$ acting on $\Z^n$ by permuting the coordinates. The lattice
$A_{n-1}^{\otimes 2}=A_{n-1}\otimes A_{n-1}$ can be replaced by lattices that 
are equivalent in a
certain specific sense --- see Sections \ref{equivalence} and \ref{formanek}
for precise formulations --- and are, hopefully, more tractable than 
$A_{n-1}^{\otimes 2}$. Indeed, most cases where the rationality problem for
generic division algebras has been successfully treated rely on such a replacement.
This is indicated in the following table, where $\sim$ denotes the aforementioned 
lattice equivalence:
\begin{center}
\renewcommand{\arraystretch}{1.25}
\begin{tabular}{|c|c|}
\hline
$n$ & $A_{n-1}^{\otimes 2}$
\\ \hline\hline
$2$ & $\sim 0$\\ \hline
$3$ & $\sim 0$ \\ \hline
$4$ & $\sim A_3^*\otimes\Z^-$ \\ \hline
$5$ or $7$ & $\sim A_{n-1}^*$ \\ \hline
\end{tabular}
\end{center}
Here, $(\,.\,)^*$ denotes dual lattices and
$\Z^-=\bigwedge^n\Z^n$ is the \emph{sign lattice} for $\S_n$.
The case $n=2$ is obvious; in fact, $A_1^{\otimes 2}\cong\Z$ is the trivial
$\S_2$-lattice. Formanek \cite{Fo}, \cite{Fo2} dealt with $n=3$ and $n=4$,
the equivalence $A_3^{\otimes 2}\sim A_3^*\otimes\Z^-$ being crucial for
$n=4$. Finally, the equivalence $A_{n-1}^{\otimes 2}\sim A_{n-1}^*$ 
is the operative fact in Bessenrodt and Le Bruyn's solution of the
rationality problem for $n=5$ and $n=7$ \cite{BL}. More recently, Beneish \cite{B} 
has given a simpler proof of this equivalence (which is known to be false for prime values
of $n>7$ \cite{BL}).

Our first contribution is in response to the ``Representation Question" of
\cite{LeB}, which asks to find $\S_n$-lattices equivalent to $A_{n-1}^{\otimes 2}$
but having lower rank. We show that, for odd values of $n$ and for $n=2$ and $n=4$,
the exterior square $\wed{A_{n-1}}$ serves that purpose, its rank being about
half that of $A_{n-1}^{\otimes 2}$. 
(For even $n\ge 6$, however, $\wed{A_{n-1}}$ is not equivalent to $A_{n-1}^{\otimes 2}$.)
The equivalence
$A_{n-1}^{\otimes 2}\sim \wed{A_{n-1}}$ covers the cases $n=2,3,4$ in the above
table, but not $n=5,7$. We do not touch here on the actual rationality problem for
the invariant field $k(\wed{A_{n-1}})^{\S_n}$ which does seem rather
formidable still. Even the case $n=4$, where $\wed{A_3}\cong A_3^*\otimes\Z^-$,
does require some effort; see \cite{Fo2} or \cite{HK}. In this connection, it is
perhaps worth remarking that the rationality problem for the invariant field of
the $\S_4$-lattice $A_3\otimes\Z^-$ appears to be open; it is in fact the only
case of a $3$-dimensional multiplicative invariant field whose rationality was left 
undecided in \cite{HK}.

As our second contribution, we confirm a conjecture of Le Bruyn \cite[p.~108]{LeB}
by showing that $n=2$ and $n=3$ are the only cases where 
$A_{n-1}^{\otimes 2}\sim 0$. 
This was previously known for prime values of $n>3$ \cite[Corollary 1(a)]{BL} and
for all $n$ that are not square-free \cite{Sa}. 
After work on the present article was completed, Le Bruyn pointed out to us that
his conjecture, expressed in terms of generic norm tori, was also proved by 
Cortella and Kunyavski\u{\i}  \cite{CK}.
The interest in faithful $\S_n$-lattices
$M\sim 0$, often called \emph{quasi-permutation} lattices, stems from the fact
that their multiplicative invariant fields $k(M)^{\S_n}$ are known to be stably
rational over $k$; see Section \ref{mult} below.
Thus, our second result is largely negative in its implications for the
rationality problem for generic division algebras. Nevertheless, it is of some
interest in its own right. In fact, the core of its proof is a nontriviality
result for the Picard groups of certain algebras of multiplicative invariants.
Such nontrivial Picard groups, far from ubiquitous, give rise to
nontrivial torsion in the Grothendieck group $G_0$  of certain polycyclic group
algebras \cite{Lo}, a phenomenon in need of further elucidation; see 
\cite{KM},\cite{Le1},\cite{Le2} for work on this topic.

In Sections \ref{equivalence} and \ref{formanek} of this article, we have collected 
the requisite background material, while Sections \ref{augmentation} and \ref{picard}
contain our main results. The latter two sections each are preceded by an introductory
paragraph laying out their contents in a more complete and technical fashion.
%===================================================

\begin{notat}
Throughout this note, $\G$ denotes a finite group, $\S_n$ is the symmetric group on
$\{1,\ldots,n\}$, and $k$ will denote a commutative field. 
Furthermore, $\otimes=\otimes_{\Z}$ and $(\,.\,)^*=\Hom_{\Z}(\,.\,,\Z)$. 
If $\H\le\G$ is a subgroup of $\G$ and $M$
a $\Z[\H]$-module then $M\ind{\H}{\G}=\Z[\G]\otimes_{\Z[\H]}M$ denotes the
induced $\Z[\G]$-module.
\end{notat}

%%%%%%%%%%%%%%%%%%%%%%%%%%%%%%%%%

\section{Equivalence of Lattices and of Fields} \label{equivalence}

%=================================
\subsection{$\G$-lattices} \label{Glattices}
We recall some standard definitions.
A \emph{$\G$-lattice} is a (left) module over the integral
group ring $\Z[\G]$ that is free of finite rank as $\Z$-module.
A $\G$-lattice $M$ is called \emph{faithful} 
if the structure map $\G\to\Aut_{\Z}(M)$ is injective.
Further, $M$ is 
\begin{itemize}
\item a \emph{permutation lattice} if $M$ has a
	$\Z$-basis, say $X$, that is permuted by $\G$; we will write such a
	lattice as $M=\Z[X]$. 
\item \emph{stably permutation} if $M\oplus P\cong Q$ holds for
	suitable permutation $\G$-lattices $P$ and $Q$.
\item \emph{permutation projective} (or \emph{invertible}) if $M$
	is a direct summand of some permutation $\G$-lattice.
\item \emph{coflasque} if $H^1(\H,M)=0$ holds for all subgroups
	$\H\le\G$ or, equivalently, $\Ext_{\Z[\G]}(P,M)=0$ holds for all
	permutation projective $\G$-lattices $P$; see \cite[Lemma 8.3]{Sw}
	for the equivalence.
\end{itemize}
The above properties are successively weaker.

\begin{lem} Let $M$, $P$, and $Q$ be $\G$-lattices, with $P$ permutation
projective and $Q$ coflasque. Suppose that there is an exact sequence
$$
0\to M\to P\to Q/\alpha Q\to 0
$$ 
with $\alpha\in\End_{\Z[\G]}(Q)$ injective.
Then there is an exact sequence of $\G$-lattices
$$
0\to M\to P\oplus Q\to Q\to 0 \ .
$$
\end{lem}

\begin{proof}
Consider the pull-back diagram
$$
\xymatrix{
 & & 0 & 0 \\
0\ar[r] & M \ar[r]\ar@{=}[d] & P \ar[r]\ar[u] 
& Q/\alpha Q \ar[r]\ar[u] & 0\\
0\ar[r] & M \ar[r] & X \ar[r]\ar[u] 
& Q \ar[r]\ar[u]& 0\\
 & & Q \ar[u]\ar@{=}[r] & Q \ar[u]_{\alpha}  \\
& & 0\ar[u] & 0\ar[u] 
}
$$
The middle column splits, by the assumptions on $P$ and $Q$, and hence 
the middle row yields the required sequence.
\end{proof}

%=================================
\subsection{Equivalence of lattices}\label{latticeequiv}
Following \cite{CTS}, two $\G$-lattices $M$ and $N$ are called 
\emph{equivalent}, written
$M\sim N$, provided there are exact sequences of $\G$-lattices
\begin{equation} \label{E:equiv}
0\to M\to E\to P\to  0
\qquad\text{and}\qquad 
0\to N\to E\to Q\to  0
\end{equation}
where $P$ and $Q$ are \emph{permutation} $\G$-lattices. For a direct proof that
this does indeed define an equivalence relation and for further background,
see \cite[Lemma 8]{CTS} or \cite{Sw}. 
Lattices $M\sim 0$ are called \emph{quasi-permutation}; explicitly,
$M$ is quasi-permutation iff there is an exact sequence of $\G$-lattices
$0\to M\to P\to Q\to  0$ where $P$ and $Q$ are permutation lattices.

In the following lemma, we collect some basic properties of the equivalence $\sim$.
All $M$, $N$ etc. are understood to be $\G$-lattices and all maps considered are
$\G$-equivariant.

\begin{lem} 
\begin{enumerate}
\item[\textnormal{(i)}]
$M\sim N\text{ and } M'\sim N'\quad\Longrightarrow\quad 
M\oplus M'\sim N\oplus N'$. 
\item[\textnormal{(ii)}]
If $0\to M\to N\to S\to 0$ is exact 
with $S$ stably permutation then $M\sim N$.
%\item[\textnormal{(iii)}]
%If $M\sim N$ and $P$ is a permutation lattice then
%$M\otimes P\sim N\otimes P$.
%\item[\textnormal{(iv)}]
%Suppose that there is an exact sequence
%$$
%0\to M\to P\to Q/\alpha Q\to 0
%$$ 
%with $\alpha\in\End_{\Z[\G]}(Q)$ injective, $P$ permutation
%projective and $Q$ stably permutation.
%Then $M\sim P$.
\end{enumerate}
\end{lem}

\begin{proof}
(i) From \eqref{E:equiv} we obtain exact sequences 
$0\to M\oplus M'\to E\oplus M'\to P\to  0$
and $0\to N\oplus M'\to E\oplus M'\to Q\to  0$; so
$M\oplus M'\sim N\oplus M'$. Similarly, $N\oplus M'\sim N\oplus N'$, and
hence $M\oplus M'\sim N\oplus N'$, by transitivity of $\sim$.

(ii) Say $S\oplus P\cong Q$ for permutation lattices $P$ and $Q$.
Then we have exact sequences 
$0\to M\to N\oplus P\to S\oplus P\cong Q\to  0$
and $0\to N\to N\oplus P\to P\to 0$, proving $M\sim N$.
%
%(iii) follows from exactness of $(\,.\,)\otimes P$ and the fact that
%tensor products of permutation modules are permutation modules.
%
%(iv) Using parts (i), (ii) and Lemma \ref{Glattices}, we obtain 
%$M\sim P\oplus Q\sim P$.
\end{proof}

%=================================
\subsection{$\G$-fields} \label{Gfields}
A \emph{$\G$-field} is a field $F$ together with a given action of $\G$ by
automorphisms on $F$. Following \cite{B}, we say that two $\G$-fields
$K$ and $F$ containing a common $\G$-subfield $k$ are 
\emph{stably isomorphic as $\G$-fields over $k$} provided, for suitable
$r$ and $s$, there is a $\G$-equivariant field isomorphism
$$
K(x_1,\ldots,x_r)\iso F(y_1,\ldots,y_s)
$$
which is the identity on $k$. Here, the $x$'s and $y$'s are transcendental
over $K$ and $F$, resp., and $K(x_1,\ldots,x_r)$ and $F(y_1,\ldots,y_s)$
are regarded as $\G$-fields with the trivial $\G$-action on the $x$'s and $y$'s.
In particular, $K(x_1,\ldots,x_r)^{\G}=K^{\G}(x_1,\ldots,x_r)$ and similarly
for $F(y_1,\ldots,y_s)$. The above isomorphism restricts to an isomorphism
$K^{\G}(x_1,\ldots,x_r)\iso F^{\G}(y_1,\ldots,y_s)$ which is the identity on 
$k^{\G}$; so:
\begin{quote}
\emph{If $K$ and $F$ are stably isomorphic as $\G$-fields over $k$ then their
fixed subfields $K^{\G}$ and $F^{\G}$ are stably isomorphic over $k^{\G}$.} 
\end{quote}

%=================================
\subsection{Multiplicative $\G$-fields} \label{mult}
We are particularly interested in $\G$-fields arising 
from $\G$-lattices $M$ 
by first extending the $\G$-action on $M$ to the group algebra $k[M]$ of $M$
over $k$, with $\G$ acting trivially on $k$, and then further to the field of fractions
$k(M)=Q(k[M])$. This particular type of $\G$-field will be called 
\emph{multiplicative}. 

More generally, the above construction can be carried out for any
$\G$-field $K$, since the actions of $\G$ on $K$ and $M$ extend uniquely to
the group algebra $K[M]$, and then to $K(M)=Q(K[M])$. If the $\G$-action on $K$
is \emph{faithful} then the field of invariants $K(M)^{\G}$ of the $\G$-field 
$K(M)$ thus obtained is called the \emph{field of tori-invariants of $K$ and
$M$ under $\G$}.

Various versions of the following proposition can be found in the literature; our
particular formulation was motivated by \cite[Lemma 2.1]{B}. 

\begin{prop}
Let $M$ and $N$ be faithful $\G$-lattices. Then 
the multiplicative $\G$-fields $k(M)$ and $k(N)$ are stably isomorphic
over $k$ if and only if $M\sim N$.
In particular, the fixed fields $k(M)^{\G}$ and $k(N)^{\G}$ are stably
isomorphic over $k$ in this case.
\end{prop}

\begin{proof}
First, $M\sim N$ implies that 
$k(M)$ and $k(N)$ are stably isomorphic as $\G$-fields over $k$.
Indeed, given the exact sequences \eqref{E:equiv}, the usual
``Galois descent" argument (e.g., \cite[Prop. 1.3]{L}) shows that
$k(E)$ is generated over $k(M)$ by algebraically independent
$G$-invariant elements, and similarly for $k(E)/k(N)$.

Conversely, if $k(M)$ and $k(N)$ are stably isomorphic as 
$\G$-fields over $k$ then, for any $\G$-field $K$ containing $k$, 
$K(M)$ and $K(N)$ are stably isomorphic as $\G$-fields over $K$. 
Indeed, for suitable $r$ and $s$, the fields $k(M\oplus{\Z}^r)$ and
$k(N\oplus{\Z}^s)$ are isomorphic as $\G$-fields over $k$, and hence 
$K(M\oplus{\Z}^r)=Q(K\otimes_kk(M\oplus{\Z}^r))$ and 
$K(N\oplus{\Z}^s)=Q(K\otimes_kk(N\oplus{\Z}^s))$ are isomorphic 
as $\G$-fields over $K$.
As noted in \ref{Gfields}, this implies that
the invariant fields $K(M)^{\G}$ and $K(N)^{\G}$ are stably isomorphic
over $K^{\G}$. In particular, choosing $K$ so that $\G$ acts faithfully
on $K$, we deduce from \cite{CTS} (or \cite[Theorem 10.1]{Sw}) that
$M\sim N$.

Since the statement about fixed fields has already been explained in
\ref{Gfields}, the proof is complete.
\end{proof}

%%%%%%%%%%%%%%%%%%%%%%%%%%%%%%%%%%%%%%%%%%%%%%%%%%%

\section{Review of the Formanek Strategy} \label{formanek}

%=================================
\subsection{Some $\S_n$-lattices} \label{lattices}
We introduce some lattices that will play an important role throughout
this note. First, $U_n$ denotes the standard permutation $\S_n$-lattice of rank $n$;
so $U_n=\Z\ind{\S_{n-1}}{\S_n}$ or, explicitly,
$$
U_n=\Z u_1\oplus\ldots\oplus\Z u_n\quad\text{with $\S_n$-action given by
$\sigma u_i=u_{\sigma(i)}$.}
$$
Sending each $u_i\mapsto 1\in\Z$, the trivial $\S_n$-lattice, yields a 
map of $\S_n$-lattices $\varepsilon: U_n\to\Z$, called the 
\emph{augmentation map}. Its kernel is the 
\emph{root lattice} $A_{n-1}$. Thus,
$$
A_{n-1}=\left\{\sum_{i=1}^n z_iu_i\Big| \sum_iz_i=0\right\}
=\Z a_1\oplus\ldots\oplus\Z a_{n-1}\quad\text{with $a_i=u_i-u_n$.}
$$
Finally, we put 
$$
A_{n-1}^{\otimes 2}=A_{n-1}\otimes A_{n-1}\ .
$$
Tensoring the exact sequence $0\to A_{n-1}\tto U_n
\stackrel{\varepsilon}{\tto}\Z\to  0$
with $A_{n-1}$ we obtain an exact sequence of $\S_n$-lattices
\begin{equation} \label{E:Gn}
0\to A_{n-1}^{\otimes 2}\tto V_n=U_n\otimes A_{n-1}\tto
A_{n-1}\to  0\ .
\end{equation}
Note that the restriction of $A_{n-1}$ to $\S_{n-1}=\stab_{\S_n}(n)$
is isomorphic to $U_{n-1}=\Z\ind{\S_{n-2}}{\S_{n-1}}$. Hence,
$V_n=\Z\ind{\S_{n-1}}{\S_n}\otimes A_{n-1}\cong\Z\ind{\S_{n-2}}{\S_n}$.

%=============================================
\subsection{The Formanek-Procesi Theorem} \label{procesi}
Let $M_n(k)$ denote the space of $n\times n$-matrices over the field $k$. The 
group $\PGLn$ operates on $M_n(k)\oplus M_n(k)$ by simultaneous conjugation.
The invariant field of this action is commonly denoted by $C_n$; so
$$
C_n=k(M_n(k)\oplus M_n(k))^{\PGLn}\ .
$$
The ultimate goal is to prove (stable) rationality of $C_n$ over $k$.
Recall that a field extension $F/k$ is called \emph{rational} if it is
purely transcendental and \emph{stably rational} if there is an extension
field $E\supseteq F$ such that $E/F$ and $E/k$ are both rational. Equivalently,
$F$ is stably rational over $k$ iff $F$ is stably isomorphic to $k$ over $k$.

The Formanek-Procesi Theorem expresses $C_n$ as the multiplicative invariant
field arising from the $\S_n$-action on the lattice $U_n\oplus U_n\oplus 
A_{n-1}^{\otimes 2}$.

\begin{thm} 
$C_n\cong k(U_n\oplus U_n\oplus A_{n-1}^{\otimes 2})^{\S_n}$
\end{thm}

This is proved in \cite[Theorem 3]{Fo} building on
Procesi's work \cite{Pr}.

%=============================================
\subsection{The strategy} \label{strategy}

Formanek's strategy to prove stable rationality of $C_n$ over $k$ consists of
finding a faithful $\S_n$-lattice $H_n$ satisfying
\begin{enumerate}
\item[(i)] $H_n\sim A_{n-1}^{\otimes 2}$, and
\item[(ii)] $k(H_n)^{\S_n}$ is stably rational over $k$.
\end{enumerate}
The above formulation of the strategy is taken from \cite{BL}. The fact
that (i) and (ii) together imply stable rationality of 
$C_n$ over $k$ is immediate from Theorem \ref{procesi} and Proposition \ref{mult},
because $U_n\sim 0$ and so 
$U_n\oplus U_n\oplus A_{n-1}^{\otimes 2}\sim A_{n-1}^{\otimes 2}$, by 
Lemma \ref{latticeequiv}(i).

%%%%%%%%%%%%%%%%%%%%%%%%%%%%%%%%%%%%%%%%%%%%%%%%%%%%%%%%%%%%%%%%%%%%%%%%%%%%%%%%

\section{The Augmentation Kernel}\label{augmentation}

%=================================
\subsection{Overview} \label{augover}
In this section, we will consider the following situation:
\begin{tabbing}
  XXXXX\=$SSSSSS$\qquad\= \kill          %sample line
  \>$U=\Z[X]$ \> will denote a permutation $\G$-lattice, \\
  \>$\varepsilon: U\to\Z$ \> is the augmentation sending all $x\in X$ to 1, and \\
  \>$A=\Ker\varepsilon$ \> will be the \emph{augmentation kernel}. 
\end{tabbing}
From the exact sequence 
\begin{equation}\label{E:augseq}
1\to A\tto U\stackrel{\varepsilon}{\tto}\Z\to 1
\end{equation}
one sees that $A$ is quasi-permutation.
It is easy to see that $A$ is coflasque or, equivalently, stably permutation
precisely if the sizes of the $\G$-orbits in $X$ are coprime.
Indeed, from the cohomology sequence that is associated with \eqref{E:augseq}
one obtains $H^1(\G,A)\cong \Z/\sum_{\O}|{\O}|\Z$, where $\O$ runs over the orbits of
$\G$ in $X$; cf. Lemma \ref{prfprop} below. 
Thus, $H^1(\G,A)=0$ iff the orbit sizes are coprime. On the other hand,
the vanishing
of $H^1(\G,A)=\Ext_{\Z[\G]}(\Z,A)$ forces the splitting of \eqref{E:augseq}; so $A$
is stably permutation in this case.
%every subgroup of $\G$
%of prime power order has a fixed point in $X$. 

By contrast, we will show
that the \emph{second symmetric power} $\sym{A}$ of $A$ is always coflasque, and is
stably permutation if there is a $\G$-orbit of odd length in $X$. As a consequence, the  
\emph{second exterior power} $\wed{A}$ of $A$ satisfies
$$
\wed{A}\sim A^{\otimes 2}=A\otimes A
$$
in this case. These facts will be proved in (\ref{stableper}) below after recalling the
requisite background concerning $\sym{A}$ and $\wed{A}$ in (\ref{symext}).

%=================================
\subsection{The second symmetric and exterior powers} \label{symext}
Let $M$ be a $\G$-lattice.
By definition, $\sym{M}$ is the quotient of 
$M^{\otimes 2}=M\otimes M$ modulo the subgroup that is generated by the
elements $m\otimes m' - m'\otimes m$ for $m,m'\in M$. We will write $mm'\in\sym{M}$ 
for the image of $m\otimes m'$; so $mm'=m'm$. A $\Z$-basis of $\sym{M}$ is given by
$\{m_im_j\mid 1\le i\le j\le r\}$, where $\{m_1,\ldots,m_r\}$
is any $\Z$-basis of $M$. Thus, 
$$
\rk\sym{M}=\binom{r+1}{2}\quad\text{with }\  r=\rk M.
$$ 
Similarly, $\wed{M}$
is the quotient of $M^{\otimes 2}$ modulo the subgroup that is generated by the
elements $m\otimes m$ for $m\in M$. Denoting the image of $m\otimes m'$ in $\wed{M}$ 
by $m\wedge m'$, a $\Z$-basis of $\wed{M}$ is given by 
$\{m_i\wedge m_j\mid 1\le i<j\le r\}$; so
$$
\rk\wed{M}=\binom{r}{2}\ .
$$ 
The action of $\G$ on $M^{\otimes 2}$ passes down to $\sym{M}$ and to $\wed{M}$,
making both $\G$-lattices. The $\G$-lattice $\wed{M}$ can be
identified with the sublattice of \emph{antisymmetric tensors} in 
$M^{\otimes 2}$, that is,
$$
\wed{M}\cong \mathsf{A}'_2(M)=\{x\in M^{\otimes 2}\mid x^{\tau}=-x\}
$$
where $\tau:M^{\otimes 2}\to M^{\otimes 2}$ is the switch 
$(m\otimes m')^{\tau}=m'\otimes m$; see
\cite[Exercise 8 on p.~A~III.190]{Bou}.
Furthermore,  $\mathsf{A}'_2(M)$ is exactly the kernel of the
canonical map $M^{\otimes 2}\onto\sym{M}$, and hence we have an exact
sequence of $\G$-lattices
\begin{equation} \label{E:exact}
0\to\wed{M}\tto M^{\otimes 2}\tto\sym{M}\to 0 \ .
\end{equation}
The foregoing, and much of the subsequent material, extends \textsl{mutatis mutandis} 
to higher tensor, symmetric, and exterior powers of permutation modules, but we will
concentrate on squares.

%=================================
\subsection{The case of permutation lattices and their augmentation kernels} 
\label{permcase} 
Recall that,
for a $\G$-set $X$, the cartesian product $X\times X$ becomes a $\G$-set with the
diagonal $\G$-action, and the set of $2$-subsets of $X$,
$$
\P(X)=\{Y\subset X\mid |Y|=2\}\ ,
$$
is canonically a $\G$-set as well. 

A $\G$-lattice $M$ is called \emph{monomial} if $M$ has a $\Z$-basis that is
permuted by $\G$ up to $\pm$-sign. Equivalently, $M$ is monomial iff $M$ is
a direct sum of lattices that are induced from rank-1 lattices for suitable
subgroups $\H\le\G$. Any $\H$-lattice of rank 1 is given by a 
homorphism $\varphi:\H\to\{\pm 1\}=\GL_1(\Z)$; we will use the notation $\Z_{\varphi}$
to denote such a lattice. Note that we have the following inclusions:
$$
\left\{\text{permutation lattices}\right\}
\subset
\left\{\text{monomial lattices}\right\}
\subset
\left\{\text{quasi-permutation lattices}\right\}\ .
$$
All these inclusions are strict; see (\ref{monomial}) below.

\begin{lem}
Let $U=\Z[X]$ be a permutation $\G$-lattice and $A$ its augmentation kernel, as in 
(\ref{augover}). Then:
\begin{enumerate}
\item[\textnormal{(i)}] 
$U^{\otimes 2}\cong
\Z[X\times X]$ and $\sym{U}\cong U\oplus \Z[\P(X)]$; all these are permutation
$\G$-lattices. Furthermore, $\wed{U}\cong\oplus_Y\Z_{\varphi_Y}\ind{\G_Y}{\G}$ 
is a monomial lattice. Here, $Y$ runs over a representative set of the
$\G$-orbits in $\P(X)$, $\G_Y=\stab_{\G}(Y)$, and $\varphi_Y:\G_Y\to\{\pm 1\}$
is the sign character for the action of $\G_Y$ on $Y$. 
\item[\textnormal{(ii)}] There is a $\G$-equivariant exact sequence
$$
0\to \sym{A}\to  \Z[\P(X)]\to A/2A\to 0\ .
$$ 
\item[\textnormal{(iii)}] There is an exact sequence of $\G$-lattices
$$
0\to \wed{A}\to  \wed{U}\to A\to 0\ .
$$ 
\end{enumerate}
\end{lem}

\begin{proof} 
(i)
The assertion about $U^{\otimes 2}$ is clear.
If $\le$ is
a fixed ordering of $X$ then $\sym{U}$ has $\Z$-basis 
$$
\{xx'\mid x\le x'\}=\{x^2\mid x\in X\}\uplus \{xx'\mid x<x'\}\cong X\uplus\P(X)
$$
as $\G$-sets. Thus, 
$$
\sym{U} = U'\oplus W \quad\text{with $U'=\sum_x\Z x^2\cong U$ and $W=\sum_{x<x'}\Z xx'
\cong\Z[\P(X)]$.}
$$
Also, $\wed{U}$ has $\Z$-basis $\{x\wedge x'\mid x<x'\}$ and $\G$ permutes this
basis up to $\pm$. The description of $\wed{U}$ follows.
%from this, and
%the fact that $\wed{U}$ is quasi-permutation follows from this description or else
%from the sequence \eqref{E:exact}.

(ii) 
The $\G$-map $X\times X\to U$, $(x,x')\mapsto x+x'$, gives rise to a map
of $\G$-lattices $U^{\otimes 2}\cong\Z[X\times X]\to U$ which passes down to
a map of $\G$-lattices 
$$
\rho: \sym{U}\tto U\ ,\qquad xx'\mapsto x+x'\quad (x,x'\in X).
$$
One easily checks that 
$\Im\rho=\varepsilon^{-1}(2\Z)$. Since $\sym{A}\subset
\sym{U}$ is spanned by the elements of the
form $(x-x')(y-y')$ for $x,x',y,y'\in X$, and 
$$
\rho((x-x')(y-y'))=x+y-(x'+y)-(x+y')+x'+y'=0 \ ,
$$ 
we have $\sym{A}\subset\Ker\rho$. Now, $\sym{A}$ is a pure sublattice of 
$\sym{U}$ and
$$
\rk\sym{U}-\rk\sym{A}
=\binom{|X|+1}{2}-\binom{|X|}{2}=|X|=\rk\Im\rho\ .
$$
This implies that, in fact, $\sym{A}=\Ker\rho$.
Thus, we obtain an exact sequence of $\G$-lattices
$$
0\to\sym{A}\stackrel{\text{incl.}}{\tto}
\sym{U}\stackrel{\rho}{\tto}
\varepsilon^{-1}(2\Z)\to 0\ .
$$
Using the above decomposition $\sym{U} = U'\oplus W$
and noting that $\rho(U')=2U$, we deduce 
the following exact sequence
$$
0\to U'+\sym{A}\stackrel{\text{incl.}}{\tto}
\sym{U}=U'\oplus W
\stackrel{\overline{\rho}}{\tto}
\varepsilon^{-1}(2\Z)/2U\to 0\ .
$$
The sum $U'+\sym{A}\subset\sym{U}$ is direct, as follows by counting ranks.
Furthermore, using the notation 
$\overline{\phantom{i}.\phantom{i}}=(\,.\,)\otimes\Z/2\Z$ 
for reduction mod $2$, we have
$\varepsilon^{-1}(2\Z)/2U=
\Ker\left(\overline{U}\stackrel{\overline{\varepsilon}}{\to}\overline{\Z}\right)
=\overline{A}$. 
Therefore, we obtain an exact sequence of $\G$-modules
\begin{equation}\label{E:Gseq}
0\to\sym{A}\tto W
\stackrel{\overline{\rho}}{\tto}\overline{A}=A/2A\to 0\ ,
\end{equation} 
proving part (ii).

(iii) Consider the map (a piece of the Koszul complex)
$$
\pi: \wed{U}\tto A\ ,\qquad x\wedge x'\mapsto x-x'\quad (x<x').
$$
This map is $\G$-equivariant and surjective, and one checks that
$\wed{A}\subset\wed{U}$ is contained in the kernel. Counting ranks,
one sees that, in fact, $\wed{A}=\Ker\pi$, which completes the proof.
\end{proof}

The sequence in part (iii) above, for $\wed{A}$, plays a role analogous to that
of sequence \eqref{E:Gn} for $A_{n-1}^{\otimes 2}$.

%=================================
\subsection{$\sym{A}$ is stably permutation} \label{stableper}
We are now ready to prove the results announced in (\ref{augover}).

\begin{prop}
The $\G$-lattice $\sym{A}$ is coflasque.
If there is a $\G$-orbit of odd length in $X$ then $\sym{A}$ is stably permutation.
Consequently, $\wed{A}\sim A^{\otimes 2}$ holds in this case.
\end{prop}

\begin{proof}
Continuing with the
notation used in the proof of Lemma \ref{permcase}, we first show that 
$\sym{A}$ is coflasque.
Inasmuch as the permutation $\G$-lattice $W$ is certainly coflasque,
the sequence \eqref{E:Gseq} reduces the assertion 
to the following statement about fixed points:
\begin{equation} \label{E:fixed}
W^{\H}\stackrel{\overline{\rho}}{\tto}\overline{A}^{\H}
\quad\text{is epi for all $\H\le\G$.}
\end{equation}
Now $\overline{A}^{\H}=\overline{U}^{\H}\cap
\Ker\left(\overline{U}\stackrel{\overline{\varepsilon}}{\to}\overline{\Z}\right)$,
and $\overline{U}^{\H}$ has $\overline{\Z}$-basis the orbit sums 
$\widehat{\O}=\sum_{x\in\O}\overline{x}$, where $\O$ is an $\H$-orbit in $X$.
Therefore, $\overline{A}^{\H}$ is spanned by the orbit sums
$\widehat{\O}$, where $\O$ has even length, together with the elements 
$\widehat{\O}+\widehat{{\O}'}$
with $\O$ and ${\O}'$ both of odd length. Both types can be written as
$\widehat{E}=\sum_{x\in E}\overline{x}$, where $E\subset X$ is an $\H$-invariant
subset of even size. For any $\H$-invariant $E\subset X$, say $E=\{x_1,\ldots,x_e\}$,
put $w_E=\sum_{1\le i<j\le e} x_ix_j$, an $\H$-invariant
element of $W$ (note that $x_ix_j=x_jx_i$). 
The image $\overline{\rho}(w_E)\in\overline{A}^{\H}$ is
$$
\sum_{1\le i<j\le e}\overline{x_i}+\overline{x_j}
= (e-1)\sum_{i=1}^e\overline{x_i}\\
= \begin{cases}
\overline{x_1}+\ldots+\overline{x_e} &\text{if $e$ is even;}\\
0 &\text{if $e$ is odd.}
\end{cases}
$$
Thus, if $e$ is even then $\overline{\rho}(w_E)=\widehat{E}$, and so
\eqref{E:fixed} is proved.

If there is an odd length $\G$-orbit $\O\subset X$ then $\overline{U}
=\overline{A}\oplus\overline{\Z}\widehat{\O}$. Thus, we can extend the map
$\overline{\rho}$ in \eqref{E:Gseq} to a $\G$-equivariant
surjection $W\oplus\Z\onto\overline{U}$ by sending $1\in\Z$ to $\widehat{\O}$.
This leads to an exact sequence of $\G$-modules
$$
0\to \sym{A}\oplus\Z\tto  W\oplus\Z
\longrightarrow \overline{U}\to 0 \ .
$$
Since $U$ and $W\oplus\Z$ are permutation lattices, Lemma \ref{Glattices} 
(with $\alpha=\cdot 2$) yields an exact sequence
\begin{equation}\label{E:newexact}
0\to \sym{A}\oplus\Z\tto  W\oplus\Z\oplus U\tto U\to 0 \ .
\end{equation}
Finally, since $\sym{A}\oplus\Z$ is coflasque, 
this sequence splits which proves that $\sym{A}$ is stably permutation.

In view of the sequence \eqref{E:exact}, the assertion about $\wed{A}$ 
now follows from Lemma \ref{latticeequiv}(ii).
\end{proof}

%=================================
\subsection{The root lattice $A_{n-1}$} \label{root}

Taking $\G=\S_n$ and $X=\{u_1,\ldots,u_n\}\cong \S_n/\S_{n-1}$, as in
(\ref{lattices}), we have
$U=U_n$ and $A=A_{n-1}$. Furthermore, $W=\bigoplus_{i<j}\Z u_iu_j\cong 
W_n:=\Z\ind{\S_{n-2}\times\S_2}{\S_n}$ and
$$
\wed{U_n}\cong \Z_{\varphi}\ind{\S_{n-2}\times\S_2}{\S_n}\ ,
$$
where $\varphi:\S_{n-2}\times\S_2\onto\{\pm 1\}$ is trivial on $\S_{n-2}$.
For odd $n$, the split sequence
\eqref{E:newexact} gives the following isomorphism of $\S_n$-lattices:
$$
\sym{A_{n-1}}\oplus U_n\oplus\Z\cong W_n\oplus U_n\oplus\Z\ .
$$
%I don't know if $\sym{A_{n-1}}$ is in fact always stably permutation. This is certainly
%so for $n=2$, because $\sym{A_1}\cong\Z$ is the trivial $\S_2$-lattice. 
For $n=2$, $\sym{A_1}\cong\Z$ is the trivial $\S_2$-lattice and, 
for $n=4$,
Formanek \cite[Lemma 5(1)]{Fo2} has shown that
$$
\sym{A_3}\oplus\Z\cong U_4\oplus\Z\ind{\mathcal{D}}{\S_4}
$$
where $\mathcal{D}$ is the Sylow 2-subgroup of $\S_4$. 

Proposition \ref{stableper} gives the equivalence
\begin{equation}\label{E:wedequiv}
\wed{A_{n-1}}\sim A_{n-1}^{\otimes 2}
\end{equation}
for odd $n$ and for $n=2,4$.
(We will see later (\ref{even}) that \eqref{E:wedequiv} fails to hold for
even $n\ge 6$.)
In view of (\ref{strategy}), the equivalence \eqref{E:wedequiv}
reduces the 
rationality problem for $C_n$, for odd $n$, to the corresponding problem for
the $\S_n$-lattice $\wed{A_{n-1}}$:
\begin{equation}\label{E:application}
\text{
\emph{Is $k(\wed{A_{n-1}})^{\S_n}$ (stably) rational over $k$?}
}
\end{equation}
Rationally, $\wed{A_{n-1}}$ is irreducible. In fact,
$$
\wed{A_{n-1}}\otimes\Q\cong S^{(n-2,1^2)} \ ,
$$
where $S^{(n-2,1^2)}$ is the \emph{Specht module} 
corresponding to the partition $(n-2,1^2)$ of $n$; e.g., \cite[Exercise 4.6 on p.~48]{FH}.
\bigskip 

Here are some small cases:

\subsubsection{The case $n=2$} Clearly, $\wed{A_1}=0$, and so 
\eqref{E:application} is trivial.

\subsubsection{The case $n=3$} Since $\wed{A_2}=\Z\cdot a_1\wedge a_2$, with $a_i$ as in 
(\ref{lattices}), one sees that
$$
\wed{A_2}\cong\Z^-\sim 0\ .
$$ 
So,
$k(\wed{A_2})=k(t)$, with $(1,2)\cdot t= t^{-1}$ and with $(1,2,3)$ acting trivially.
The invariant field is $k(t+t^{-1})$, a rational extension of $k$.
(This is also clear from L\"uroth's Theorem.)

\subsubsection{The case $n=4$} Here,
$$
\wed{A_3}\cong\Hom_{\Z}(A_3,\Z^-)\ .
$$
This isomorphism is a consequence of the following equalities of $\S_4$-lattices:
$$
\left(\wed{A_3}\right)\wedge A_3=\bigwedge^3 A_3=\Z_{\det A_3}= \Z^- \ .
$$
Further, $\Hom_{\Z}(A_3,\Z^-)\cong A_3^*\otimes\Z^-\cong (SA_3)^*$, 
the dual of the \emph{signed root lattice}
$SA_3=\Z^-\otimes A_3$. Formanek has shown \cite[Theorems 13 and 14]{Fo2}
that $k(SA_3^*)^{\S_4}$ is rational over $k$. The same result is also
covered as case $W_8(198)$ in \cite{HK}.

%%%%%%%%%%%%%%%%%%%%%%%%%%%%%%%%%%%%%%%%%%%%%%%%%%%

\section{Picard Groups} \label{picard}

%=================================
\subsection{Statement of results} \label{statement}
Following \cite{CTS,CTS2} we put, for any $\Z[\G]$-module $M$, 
$$
\cyr{Sh}^i(\G,M)=\bigcap_{g\in\G} 
\Ker\left(\Res^{\G}_{\gen{g}}:H^i(\G,M)\to H^i(\gen{g},M)\right)\ .
$$
Of particular interest for us will be the case where $M$ is a $\G$-lattice
and $i=1$ or $2$. Indeed, if $\sim$ denotes the equivalence of (\ref{latticeequiv}) then
$$
\cyr{Sh}^2(\G,\,.\,)\ \text{\emph{is constant on $\sim$-classes;}}
$$
see Lemma~\ref{connection} below. Furthermore,
by \cite{Lo}, one has
$$
\cyr{Sh}^1(\G,M)=\Pic(k[M]^{\G})\ ,
$$
the \emph{Picard group} of the algebra of invariants under the (multiplicative)
action of $\G$ on the group algebra $k[M]$. These Picard groups are seldom
nontrivial: If $\rk M=3$ then only 2 out of the 73 conjugacy classes of finite
subgroups $\G\le\GL_3(\Z)$ lead to invariants with nontrivial Picard group, and
only 9 out of 710 occur for $\rk M=4$; see \cite{Lo}. Nevertheless, we
will prove here the following nontriviality result for the root lattice $M=A_{n-1}$.
Note that, since $A_{n-1}\sim 0$, the result shows in particular that
$\cyr{Sh}^1(\G,\,.\,)$ is not in general a $\sim$-invariant.

\begin{prop}
\begin{enumerate}
\item[\textnormal{(i)}]
Suppose $n$ is not prime and let $p$ be the smallest prime divisor of $n$. Then there 
exists a subgroup $\H\cong C_p\times C_p$ of
$\S_n$ such that $\cyr{Sh}^1(\H,A_{n-1})\cong\Z/p\Z$.
\item[\textnormal{(ii)}]
If $n=2^a\ge 4$ then $\cyr{Sh}^1(\A_n,A_{n-1})\cong\cyr{Sh}^1(\T,A_{n-1})
\cong\Z/2\Z$, where $\A_n$ is the alternating group of degree $n$ and $\T$
denotes the Sylow 2-subgroup of $\A_n$.
\end{enumerate}
\end{prop}

\noindent Part (ii), with $n=4$, yields the two nontrivial Picard groups in rank 3.
The following corollary of the proposition confirms a conjecture of
Le Bruyn \cite[p. 108]{LeB}. 

\begin{cor}
If $n$ is not prime then the $\S_n$-lattice $A_{n-1}^{\otimes 2}$ 
is not equivalent to a direct
summand of a quasi-permutation lattice.
\end{cor}

\noindent This corollary will be immediate from part (i) of the Proposition together with
Lemma~\ref{connection} below. Inasmuch as $A_{n-1}^{\otimes 2}$ is known
to be permutation projective for prime values of $n$ \cite[Proposition 3]{BL},
our non-primeness hypothesis on $n$ is definitely necessary. We remark that
$A_{n-1}^{\otimes 2}$ was previously known not to be quasi-permutation for all $n$ that
are divisible by a square \cite{Sa} or whenever
$n$ is a prime $>3$ \cite[Corollary 1(a)]{BL}.
Summarizing, we now know:
\medskip

\centerline{
\emph{The $\S_n$-lattice $A_{n-1}^{\otimes 2}$ is quasi-permutation precisely for $n=2$
and $n=3$.}
}
\medskip

\noindent Finally, we consider even values of $n\ge 6$ and show in (\ref{even})
that $\wed{A_{n-1}}$ is no longer equivalent to $A_{n-1}^{\otimes 2}$ for such $n$.

%=================================
\subsection{Connection with equivalence of lattices} \label{connection}
Every $\G$-lattice $M$ has a \emph{flasque resolution}
\begin{equation} \label{E:flasqueres}
0\to M\to P\to F\to  0
\end{equation}
with $P$ a permutation $\G$-lattice and $F$ a \emph{flasque} $\G$-lattice,
that is, $H^{-1}(\H,F)=0$ holds for all subgroups $\H\le\G$. Moreover, $F$ is
determined by $M$ up to stable equivalence: If
$0\to M\to P'\to F'\to  0$
is another flasque resolution of $M$ then $F$ and $F'$ are \emph{stably equivalent},
that is, $F\oplus Q\cong F'\oplus Q'$ holds for suitable permutation $\G$-lattices
$Q$ and $Q'$. Following \cite{CTS}, the stable equivalence class of $M$ will be 
written $[M]$, and the stable equivalence class of $F$ in the flasque resolution
\eqref{E:flasqueres} will be denoted $\rho(M)$; so
$$
\rho(M)=[F]\ .
$$
Finally, for any two $\G$-lattices $M$ and $N$,
$$
M\sim N \iff \rho(M)=\rho(N)\ .
$$
See \cite[Lemme 5]{CTS} for all this.

Note that $H^{\pm 1}(\G,M)$ depends only on the stable equivalence class $[M]$, because
$H^{\pm 1}$ is trivial for permutation modules. The following lemma is extracted from
\cite[pp. 199--202]{CTS2}. 

\begin{lem}
\begin{enumerate}
\item[\textnormal{(i)}] Let $0 \to M\to P\to N\to 0$ be an exact sequence of $\Z[\G]$-modules,
with $P$ a permutation projective $\G$-lattice. Then $\cyr{Sh}^2(\G,M)\cong\cyr{Sh}^1(\G,N)$.
\item[\textnormal{(ii)}] For any $\G$-lattice $M$, $H^1(\G,\rho(M))\cong\cyr{Sh}^2(\G,M)$.
\item[\textnormal{(iii)}] If $M$ is equivalent to a direct summand of a quasi-permutation 
$\G$-lattice then $\cyr{Sh}^2(\H,M)=0$ holds for all subgroups $\H\le\G$. 
\end{enumerate}
\end{lem}

\begin{proof}
(i) First, $\cyr{Sh}^2(\G,P)=0$, because
$\cyr{Sh}^i(\G,\,.\,)$ is additive on direct sums and $H^2(\G,\Z\ind{\H}{\G})
\cong\Hom(\H,\Q/\Z)$ is detected by restrictions to cyclic subgroups.
The isomorphism $\cyr{Sh}^2(\G,M)\cong\cyr{Sh}^1(\G,N)$ now follows from the 
commutative exact diagram below, obtained from the cohomology
sequences that are associated with the given exact sequence:
$$
\xymatrix{
0=H^1(\G,P)\ar[r] & H^1(\G,N) \ar[r]\ar[d]^{\Res} & H^2(\G,M) \ar[r]\ar[d]^{\Res}  
& H^2(\G,P) \ar@{ >->}[d]^{\Res} \\
0=\prod_g H^1(\gen{g},P)\ar[r] & \prod_g H^1(\gen{g},N)\ar[r] & \prod_g H^2(\gen{g},M) \ar[r] & 
\prod_g H^2(\gen{g},P)
}
$$

(ii) Let $0\to M\to P\to F\to  0$ be a flasque resolution of $M$; so $\rho(M)=[F]$.
By periodicity of cohomology for cyclic groups, we have
$H^1(\gen{g},\rho(M))\cong H^{-1}(\gen{g},\rho(M))=0$
for all $g\in\G$, because $F$ is flasque. Thus, $\cyr{Sh}^1(\G,\rho(M))=H^1(\G,\rho(M))$
and the isomorphism $H^1(\G,\rho(M))\cong\cyr{Sh}^2(\G,M)$ follows from part (i). 

(iii) Since the property of being equivalent to a direct summand of a quasi-permutation 
lattice is inherited by restrictions
to subgroups, it suffices to show that $\cyr{Sh}^2(\G,M)=0$ 
or, equivalently, $H^1(\G,\rho(M))=0$ holds whenever $M\oplus M'$ is quasi-permutation.
But $\rho(\,.\,)$ is easily seen to be additive, as is $H^1(\G,\,.\,)$. So
$H^1(\G,\rho(M))\oplus H^1(\G,\rho(M'))\cong H^1(\G,\rho(M\oplus M'))=0$,
because $\rho(M\oplus M')=[0]$. This entails that $H^1(\G,\rho(M))=0$, as desired.
\end{proof}

\textit{Deduction of Corollary \ref{statement} from Proposition \ref{statement}(i).} 
Applying part (i) of the Lemma to the sequence \eqref{E:Gn}, we deduce the isomorphism
\begin{equation} \label{E:iso}
\cyr{Sh}^2(\H,A_{n-1}^{\otimes 2})\cong\cyr{Sh}^1(\H,A_{n-1})
\end{equation}
for any subgroup $\H\le\S_n$. Inasmuch as $\cyr{Sh}^1(\H,A_{n-1})$ is nontrivial 
for certain $\H$
if $n$ is not prime, by Proposition \ref{statement}(i), we conclude from part (iii) of
the Lemma that the $\S_n$-lattice $A_{n-1}^{\otimes 2}$ is not equivalent to a direct
summand of a quasi-permutation lattice, as asserted in Corollary \ref{statement}.
\hfill$\square$

%=================================
\subsection{Proof of Proposition \ref{statement}} \label{prfprop}
We begin with a simple lemma which shows in particular that the restriction
map $H^1(\S_n,A_{n-1})\to H^1(\gen{\tau},A_{n-1})$ is an isomorphism for
$\tau=(1,2,\ldots,n)$; so $\cyr{Sh}^1(\S_n,A_{n-1})=\{0\}$.

\begin{lem} For any subgroup $\H\le\S_n$, 
$H^1(\H,A_{n-1})\cong \Z/\sum_{\O}|{\O}|\Z$, where $\O$ runs over the orbits of
$\H$ in $\{1,\ldots,n\}$. Moreover, the restriction map
$H^1(\S_n,A_{n-1})\to H^1(\H,A_{n-1})$ is surjective.
\end{lem}

\begin{proof}
From the cohomology sequence that is associated with the augmentation
sequence $1\to A_{n-1}\to U_n=\Z[\{1,\ldots,n\}]\stackrel{\varepsilon}{\to}\Z\to 1$,
one obtains the exact sequence
$U_n^{\H}\stackrel{\varepsilon}{\to}\Z\to H^1(\H,A_{n-1})\to 0$
which implies the asserted description of $H^1(\H,A_{n-1})$.
For surjectivity of restriction maps, consider the commutative diagram
of cohomology sequences
$$
\xymatrix{
U_n^{\S_n} \ar[r]\ar@{_{(}->}[d] & \Z \ar[r]\ar@{=}[d] 
& H^1(\S_n,A_{n-1}) \ar[d]^{\Res} \ar[r] & 0 \\
U_n^{\H} \ar[r] & \Z \ar[r] &  H^1(\H,A_{n-1}) \ar[r] & 0 
}
$$

\end{proof}

\textit{Proof of Proposition \ref{statement}(i).} 
Let $p$ be the smallest prime divisor of $n$; so
$n=pk$ with $k\ge p$. Define commuting $p$-cycles 
$\sigma_1,\ldots,\sigma_k\in\S_n$ by
$$
\sigma_i=((i-1)p+1,(i-1)p+2,\ldots,ip)\ .
$$
First assume $k\ge p+1$. Define $\alpha,\beta\in\S_n$ by
$$
\alpha=\prod_{i=1}^{k-1}\sigma_i
\qquad
\text{and}
\qquad
\beta=\prod_{i=1}^{p-1}\sigma_i^{-i}\cdot
\prod_{i=p+1}^k\sigma_i \ .
$$
Put $\H=\gen{\alpha,\beta}\le\S_n$, an elementary abelian $p$-group of rank 2.
The $\H$-orbits in $\{1,2,\ldots,n\}$ are exactly the subsets 
$\O_i=\{(i-1)p+1,(i-1)p+2,\ldots,ip\}$, all of length $p$. Therefore, by the Lemma,
$$
H^1(\H,A_{n-1})\cong\Z/p\Z \ .
$$
The cyclic subgroups $\mathcal{C}\le\H$ are: $\gen{1}$, $\gen{\alpha}$ (acts trivially
on $\O_k$), and $\gen{\alpha^i\beta}\ (1\le i\le p)$ (acts trivially on $\O_i$).
Since they all have fixed points, $H^1(\mathcal{C},A_{n-1})$ is trivial in each case.
Therefore, $\cyr{Sh}^1(\H,A_{n-1})=H^1(\H,A_{n-1})=\Z/p\Z$, as required.

Now assume $n=p^2$ and define commuting $p$-cycles $\pi_1,\ldots,\pi_p\in\S_n$ by
$$
\pi_i=(i,i+p,\ldots,i+p(p-1))\ .
$$
(Imagining the numbers $1,\ldots,p^2$ arranged in a square with $i$-th row the
entries of $\sigma_i$, the $\pi_i$ are defined columnwise.) 
Put $\pi=\prod_{i=1}^p\pi_i$.
Then $\pi\sigma_i\pi^{-1}=\sigma_{i+1}$ (indices $\mod p$), and so $\pi$
commutes with $\sigma=\prod_{i=1}^p\sigma_i$. Put $\H=\gen{\sigma,\pi}\le\S_n$, 
so $\H\cong C_p\times C_p$. Note that $\H$ acts transitively on $\{1,\ldots,p^2\}$,
and hence $\{1,\ldots,p^2\}\cong\H$ as left $\H$-sets. Therefore, the orbits of
every nonidentity element of $\H$ all have length $p$. The Lemma now gives
$$
H^1(\H,A_{n-1})\cong\Z/p^2\Z
\qquad\text{and}\qquad
H^1(\gen{g},A_{n-1})\cong\Z/p\Z
$$
for each nonidentity $g\in\H$. Since the restriction maps $H^1(\H,A_{n-1})\to
H^1(\gen{g},A_{n-1})$ are surjective, by the Lemma, the kernel is the socle of
$\Z/p^2\Z$ in each case. This socle is therefore equal to $\cyr{Sh}^1(\H,A_{n-1})$,
which completes the proof of part (i).\bigskip

\textit{Proof of Proposition \ref{statement}(ii).} 
Now assume $n=2^a\ge 4$. The Lemma implies that $H^1(\A_n,A_{n-1})\cong\Z/n\Z$ is a 
$2$-group. Thus, the restriction map $H^1(\A_n,A_{n-1})\to H^1(\T,A_{n-1})$
is an isomorphism for the Sylow 2-subgroup $\T$ of $\A_n$; cf. \cite[p. 84]{Br}.
Therefore, $\cyr{Sh}^1(\A_n,A_{n-1})$ embeds into $\cyr{Sh}^1(\T,A_{n-1})$.
We may assume that $\tau=(1,2,\ldots,n)^2\in\T$. Since $H^1(\T,A_{n-1})\cong\Z/n\Z$
maps onto $H^1(\gen{\tau},A_{n-1})\cong \Z/\frac{n}{2}\Z$, we see that
$\cyr{Sh}^1(\T,A_{n-1})$ has order at most $2$. On the other hand, since $\A_n$
contains no $n$-cycle, all $H^1(\gen{g},A_{n-1})$ with $g\in\A_n$ are cyclic
$2$-groups of order smaller than $n$. Hence the socle of $H^1(\A_n,A_{n-1})\cong\Z/n\Z$ 
maps to $0$ under all restrictions $H^1(\A_n,A_{n-1})\to H^1(\gen{g},A_{n-1})$, and
so the socle is contained in $\cyr{Sh}^1(\A_n,A_{n-1})$. Therefore,
$\cyr{Sh}^1(\A_n,A_{n-1})$ and $\cyr{Sh}^1(\T,A_{n-1})$ both have order $2$, and
the proof of the Proposition is complete. \hfill$\square$

%=================================
\subsection{Monomial lattices} \label{monomial}
This section contains some observations about monomial lattices; see (\ref{permcase}).
In particular, we will show that their $H^1$ is always an elementary abelian $2$-group.
Therefore, the quasi-permutation $\S_n$-lattice $A_{n-1}$ is not
monomial for $n>2$.

\begin{lem}
Let $M$ be a monomial $\G$-lattice; so $M$ is a direct sum of lattices
isomorphic to $\Z_{\varphi}\ind{\H}{\G}$ for $\H\le\G$ and 
$\varphi:\H\to\{\pm 1\}$. Then $H^1(\G,M)\cong (\Z/2\Z)^m$, where $m$
is the number of nontrivial $\varphi$'s occuring in the sum. Furthermore,
$\cyr{Sh}^1(\G,M)=\{0\}$ and $\cyr{Sh}^2(\G,M)=\{0\}$.
\end{lem}

\begin{proof}
Since monomial lattices are quasi-permutation, we know by
Lemma~\ref{connection}(iii) that $\cyr{Sh}^2(\G,M)=\{0\}$; so we may concentrate on
$H^1$ and $\cyr{Sh}^1$.
As both functors are additive on direct sums, we may assume that
$M=\Z_{\varphi}\ind{\H}{\G}$ with $\varphi:\H\to\{\pm 1\}$ nontrivial (otherwise,
$M$ is a permutation module and the assertions are clear). Then
$$
H^1(\G,M)\cong H^1(\H,\Z_{\varphi})\cong H^1(\H/\Ker\varphi,\Z_{\varphi})=\Z/2\Z\ ,
$$
where the first isomorphism comes from Shapiro's Lemma and the second is given
by inflation; see \cite[34.2 and 35.3]{Ba}. This proves the asserted
isomorphism for $H^1(\G,M)$.

The Shapiro isomorphism 
$H^1(\G,M)\stackrel{\cong}{\tto} H^1(\H,\Z_{\varphi})$ is equal to the composite
$$
\operatorname{proj}^*\circ\Res_{\H}^{\G}:H^1(\G,M)\to H^1(\H,M)\to
H^1(\H,\Z_{\varphi})\ ,
$$
where $\operatorname{proj}:M=\Z_{\varphi}\ind{\H}{\G}\onto
\Z_{\varphi}$ is the projection onto the direct summand $1\otimes\Z_{\varphi}\cong
\Z_{\varphi}$ of $M$. Fixing $g\in\H$ with $\varphi(g)=-1$, the restriction map
$H^1(\H,\Z_{\varphi})=\Z/2\Z\to  H^1(\gen{g},\Z_{\varphi})=\Z/2\Z$ is an
isomorphism, and hence so is the map 
$\Res^{\H}_{\gen{g}}\circ\operatorname{proj}^*\circ\Res_{\H}^{\G}
=\operatorname{proj}^*\circ\Res^{\G}_{\gen{g}}$. This proves that
$\Res^{\G}_{\gen{g}}:H^1(\G,M)\to H^1(\gen{g},M)$ is injective, whence 
$\cyr{Sh}^1(\G,M)=\{0\}$. 
\end{proof}

%=================================
\subsection{The case of even $n\ge 6$} \label{even}
As promised in (\ref{root}), we now use the above techniques to prove
the following

\begin{lem}
For even $n\ge 6$, the $\S_n$-lattice $\sym{A_{n-1}}$ is not equivalent to
a direct summand of a quasi-permutation lattice.
Moreover, $\wed{A_{n-1}}\nsim A_{n-1}^{\otimes 2}$.
\end{lem}

\begin{proof}
We first deal with $\sym{A_{n-1}}$.
By Lemma~\ref{connection}(iii), it suffices to show that
$\cyr{Sh}^2(\H,\sym{A_{n-1}})\neq 0$ holds for some subgroup $\H\le\S_n$. 
To this end, we use the exact sequence of Lemma~\ref{permcase}(ii)
which now takes the following form (cf.~(\ref{root}))
$$
0\to\sym{A_{n-1}}\tto \Z\ind{\S_{n-2}\times\S_2}{\S_n} 
\tto\overline{A_{n-1}}=A_{n-1}/2A_{n-1}\to 0\ ,
$$ 
where $\overline{\phantom{i}.\phantom{i}}=(\,.\,)\otimes\Z/2\Z$ 
denotes reduction mod $2$. By Lemma \ref{connection}(i), this sequence
implies $\cyr{Sh}^2(\H,\sym{A_{n-1}})\cong \cyr{Sh}^1(\H,\overline{A_{n-1}})$, and
hence the issue is to show that 
$\cyr{Sh}^1(\H,\overline{A_{n-1}})\neq 0$ holds for some subgroup $\H\le\S_n$. 
Now, suppose we can find a subgroup $\H$ satisfying
\begin{enumerate}
\item $\H$ is a $2$-group;
\item $\H$ has no fixed point in $\{1,2,\ldots,n\}$, but
  every $g\in\H$ does have a fixed point.
\end{enumerate}
Then the cohomology sequence that is associated with the
exact sequence $0\to \overline{A_{n-1}}\to\overline{U_n}
\stackrel{\overline{\varepsilon}}{\to}\overline{\Z}\to 0$ 
yields the following commutative diagram
$$
\xymatrix{
\overline{U_n}^{\H} \ar[r]^{0}\ar[d]^{\Res} & \overline{\Z} \ar@{ >->}[r]\ar[d]^{\Res}  
& H^1(\H,\overline{A_{n-1}}) \ar[d]^{\Res} \\
\prod_g \overline{U_n}^{\gen{g}}\ar@{->>}[r] & \prod_g \overline{\Z} \ar[r]^-{0}  
& \prod_g H^1(\gen{g},\overline{A_{n-1}}) 
}
$$
Here, the first $0$ results from the fact that $\overline{U_n}^{\H}$ is spanned by the
orbit sums of the $\H$-orbits in $\{1,\ldots,n\}$, all of which have positive even
length, while the second $0$ is a consequence of the existence of $g$-fixed points, 
which implies surjectivity of 
$\overline{U_n}^{\gen{g}}\stackrel{\varepsilon}{\to}\overline{\Z}$ for all $g\in\H$.
Thus, we conclude that the (nonzero) image of $\overline{\Z}$ in $H^1(\H,\overline{A_{n-1}})$
does indeed belong to $\cyr{Sh}^1(\H,\overline{A_{n-1}})$, forcing this group to be
non-trivial as required.

The actual construction of a suitable subgroup $\H$ is a simple matter. For
example, 
we can take the group $\H=\gen{\alpha,\beta}$, with  $\alpha=(1,2)(3,4)\cdots(n-3,n-2),
\beta=(1,2)(5,6)\cdots(n-1,n)$, that was used in the proof of
Proposition ~\ref{statement}(i) for $p=2, k\ge 3$; see (\ref{prfprop}).
%putting $\alpha=(1,2)(3,4)\cdots(n-3,n-2),
%\beta=(3,4)(5,6)\cdots(n-1,n)\in\S_n$, we can take $\H=\gen{\alpha,\beta}$.
This proves the assertion about $\sym{A_{n-1}}$.

In order to show that $\wed{A_{n-1}}\nsim A_{n-1}^{\otimes 2}$, it suffices to
exhibit a subgroup $\H\le\S_n$ such that $\cyr{Sh}^2(\H,\wed{A_{n-1}})\ncong
\cyr{Sh}^2(\H,A_{n-1}^{\otimes 2})$ In fact, the above $\H$ will do. For, by
\eqref{E:iso} and Proposition~\ref{statement}(i) (see also (\ref{prfprop})), 
we know that 
$\cyr{Sh}^2(\H,A_{n-1}^{\otimes 2})\cong\cyr{Sh}^1(\H,A_{n-1})=\Z/2\Z$. 
On the other hand, 
$\H\le\stab_{\S_n}(\{n-1,n\})=\S_{n-2}\times\S_2$.
Thus, it suffices to show that
$$
\cyr{Sh}^2(\H,\wed{A_{n-1}})=\{0\}
$$
holds for every $\H\le\S_{n-2}\times\S_2$. Consider the exact sequence
$$
0\to \wed{A_{n-1}}\to  \wed{U_n}\to A_{n-1}\to 0 
$$
from Lemma~\ref{permcase}(iii) and the resulting map $\cyr{Sh}^2(\H,\wed{A_{n-1}})\to
\cyr{Sh}^2(\H,\wed{U_n})=\{0\}$. Here, the last equality holds since $\wed{U_n}$ is
monomial; see Lemmas \ref{permcase}(i) and \ref{monomial}. Thus, in order to show
that $\cyr{Sh}^2(\H,\wed{A_{n-1}})=\{0\}$, it suffices to show that
$H^2(\H,\wed{A_{n-1}})\to H^2(\H,\wed{U_n})$ is mono or, equivalently,
that
$$
H^1(\H,\wed{U_n})\to H^1(\H,A_{n-1}) \text{ is epi.}
$$
But, as $\S_{n-2}\times\S_2$-lattices (and hence as $\H$-lattices), 
$\wed{U_n}\cong \Z_{\varphi}\ind{\S_{n-2}\times\S_2}{\S_n}$ contains 
$\Z\cdot u_n\wedge u_{n-1}$
as a direct summand, and $A_{n-1}$ contains $\Z\cdot(u_n-u_{n-1})$ as a sublattice, 
with 
$$
A_{n-1}/\Z\cdot(u_n-u_{n-1})\cong\Z\ind{\S_2}{\S_{n-2}\times\S_2} \text{ a permutation
lattice.}
$$ 
Moreover, $\Z\cdot u_n\wedge u_{n-1}\stackrel{\cong}{\tto}\Z\cdot(u_n-u_{n-1})$
under $\wed{U_n}\to A_{n-1}$. Therefore, restricting $H^1(\H,\wed{U_n})\to H^1(\H,A_{n-1})$
to the subgroup $H^1(\H,\Z\cdot u_n\wedge u_{n-1})$, we obtain
$H^1(\H,\Z\cdot u_n\wedge u_{n-1})\stackrel{\cong}{\tto} H^1(\H,\Z\cdot(u_n-u_{n-1}))
\onto H^1(\H,A_{n-1})$. This completes the proof.
\end{proof}

%%%%%%%%%%%%%%%%%%%%%%%%%%%%%%%%%%%%%%%%%%%%%%%%%%%%%%%%%%%%%%%%%%%%%%%%%%%%%%%%

%\begin{ack} The authors would like to thank 
%\end{ack}

%%%%%%%%%%%%%%%%%%%%%%%%%%%%%%%%%%%%%%%%%%%%%%%%%%%%%%%%%%%%%%%%%%%%%%%%%%%%%%%%

\end{document}